\documentclass[11pt]{amsart}
\usepackage{amssymb,amscd, amsmath}

\newcommand{\cv}{{\mathcal V}}
\newcommand{\cw}{{\mathcal W}}
\newcommand{\cz}{{\mathcal Z}}
\newcommand{\fB}{{\mathfrak B}}
\newcommand{\fT}{{\mathfrak T}}
\newcommand{\fU}{{\mathfrak U}}
\newcommand{\fb}{{\mathfrak b}}
\newcommand{\fg}{{\mathfrak g}}
\newcommand{\fh}{{\mathfrak h}}
\newcommand{\fH}{{\mathfrak H}}

\newcommand{\fe}{{\mathfrak e}}

\newcommand{\mc}{{\mathbb C}}

\newcommand{\mP}{{\mathbb P}}

\newcommand{\p}{{\p}}

\newcommand{\la}{{\lambda}}

\newcommand{\A}{{\alpha}}

\newcommand{\z}{{\zeta}}

\def\mch{{\mathcal H}}
\def\ch{\mch_\fB^*}
\def\tec{H_\fT^*}
\def\ec{equivariant cohomology}
\def\coh{cohomology}
\def\eq{equivariant}
\def\rv{$\fB$-regular variety}
\def\rvs{$\fB$-regular varieties}
\def\iff{if and only if }
\def\bi{$\fB$-invariant}

\def\bisv{$\fB$-invariant subvariety}
\def\bisvs{$\fB$-invariant subvarieties}

\def\st{such that }
\def\ni{\noindent}
\def\lie{\textup{Lie}}

\def\pr{principal}
\def\psv{principal subvariety}

\def\p{\varphi}

\def\Tr{\textup{Tr}}
\def\End{\textup{End}}

\theoremstyle{plain}
\newtheorem{theorem}{Theorem}
\newtheorem{lemma}[theorem]{Lemma}
\newtheorem{corollary}[theorem]{Corollary}

\theoremstyle{definition}

\theoremstyle{remark}

\begin{document}
%\begin{center}
%\bigskip

%{\Large\bf On the \eq\ cohomology of subvarieties\\
%of a $\fB$-regular variety}

%\bigskip{\sc James B.\ Carrell and Kiumars Kaveh}

%\end{center}

\title{On the equivariant cohomology of subvarieties of a $\fB$-regular variety}
\author{\sc James B.\ Carrell and Kiumars Kaveh}
\maketitle

\medskip
\begin{center} {\small To Bert Kostant on his 80th birthday.}
\end{center}

\bigskip
\baselineskip10pt

%{\bf Abstract}
\noindent
{\tiny {\bf Abstract.}
By a \rv, we mean a smooth
projective variety over $\mc$ admitting an algebraic action of the
upper triangular Borel subgroup  $\fB
\subset SL_2(\mc)$ \st  the unipotent radical in $\fB$
has a unique fixed point. A result of M. Brion and the first
author
\cite{BC} describes the equivariant \coh\ algebra (over $\mc$)
of a \rv $X$ as the coordinate ring
of a remarkable affine curve in $X\times \mP^1$.
The main result of this paper uses this fact to
classify the $\fB$-invariant subvarieties $Y$ of a \rv
$X$ for which the restriction map $i_Y:H^*(X)\to H^*(Y)$
is surjective.}
\baselineskip12pt

\section{Introduction}
\bigskip
A nilpotent element $e$ in the Lie algebra $\fg$
of a complex semi-simple Lie group $G$ is
{\em regular} if it lies in a unique Borel subalgebra
$\fb$ of $\fg$. If we let $B$ be the unique Borel
subgroup of $G$ with $\lie(B)=\fb$
and recall that the flag variety $G/B$ of $G$
parameterizes the family of all
Borel subalgebras of $\fg$,
it follows that the one parameter group $\exp(te)$ ($t\in \mc$) of $G$
acts on $G/B$ by left translation with
unique fixed point the identity coset $B$ or, equivalently,
the unique Borel subalgebra $\fb$ containing $e$.
By the Jacobson-Morosov Lemma,
$e$ determines a two dimensional solvable subalgebra $\fe$ of $\fg$
isomorphic to the Lie algebra of the
upper triangular Borel  subgroup $\fB$ of $SL_2(\mc)$.
The two dimensional solvable subgroup $E$ of $G$
determined by $\fe$ is thus isomorphic to a
Borel either in $SL_2(\mc)$
or $\mP SL_2(\mc)$ and hence is a
homomorphic image of $\fB$. Consequently, $\fB$
acts on $G/B$, via $E$, such that its unipotent
radical $\fU$ has exactly one fixed point,
namely $B$.

In view of this, one may  generalize the notion of a regular nilpotent
by considering an algebraic action $\fB \circlearrowright X$
of the upper triangular group $\fB$
on a smooth complex projective variety $X$ \st
its unipotent radical $\fU$ has a
unique fixed point $o\in X$. For such an action,
the maximal torus $\fT$ on the diagonal of $\fB$
is known to have a finite fixed point set (see \cite{CRELLE}).
In \cite{BC,CRELLE}, such an action is called {\em regular}
and $X$ is called a regular variety.
Slightly changing this terminology,
we will henceforth call the action $\fB$-regular
and say that $X$ is a {\em \rv}. By the above remarks,
the flag variety $G/B$ and, more generally, all algebraic
homogeneous spaces $G/P$, $P$
a parabolic in $G$, give a rich class of
\rvs.

The main goal of this paper is to study the cohomology algebras
of \bisvs\ of a \rv.
For example, Schubert varieties  in a $G/B$ or $G/P$
(that is, closures of $B$-orbits)
form an important class of
examples of \bi\ subvarieties to which our results will apply.
%In Section 4, we display a class of $\fB$-regular subvarieties
%of $G/B$ which contains all the smooth Schubert varieties
%and make a remark about the relation of their Poincar\'e polynomials
%to a formula of Kostant and Macdonald.
Moreover, certain nilpotent Hessenberg
varieties  (including the Peterson variety) \cite{T}
and certain Springer fibres corresponding to nilpotents in the
centralizer of a regular nilpotent $e\in \fg$ give further
interesting examples of \bi\ subvarieties in $G/B$
which we hope to investigate in a future work.

If $X$ is a \rv, then one knows the remarkable fact
that its \coh\
algebra $H^*(X)$ over $\mc$
is isomorphic with the coordinate ring $A(Z)$
of the zero scheme $Z$ associated to the algebraic vector field on $X$
generated by $\fU$ (cf. \cite{AC2,AC,CRELLE}). Moreover,
its $\fT$-equivariant cohomology $H^*_\fT(X)$ is
isomorphic to the coordinate ring
of a canonical $\fT$-stable
affine curve $\cz_X$ in $X\times \mP^1$ (\cite{BC}).
That is, Spec$\big(H^*_\fT(X)\big)$ is $\cz _X$.
These isomorphisms extend to
what we will call {\em principal subvarieties}:
namely, \bi\ subvarieties $Y$
of $X$ for which the natural restriction map
$i_Y^*:H^*(X)\to H^*(Y)$ is surjective (cf. Theorem 1).
In particular, if $Y$ is
principal, then Spec$\big(H^*_\fT(Y)\big)$ is
the union  $\cz_Y$ of the irreducible components of
$\cz_X$ lying over $Y\times \mc$.
Furthermore, if $Y\cap Z$ denotes the schematic intersection of $Y$ and $Z$,
then the coordinate ring $A(Y\cap Z)$ is isomorphic to $H^*(Y)$ as long
as $\dim A(Y\cap Z)=\dim H^*(Y)$.
With these facts as motivation, our aim  is to
classify the principal subvarieties of a \rv.
This will follow from a general result which
describes the image of $\tec(X)$ in $\tec(Y)$
under $i_Y^*$ for any \bisv\ $Y$ of a \rv\ $X$.

The curve $\cz_X$ admits a natural description.
Consider  the diagonal action of $\fB$
on $X\times \mP^1$, where $\fB$ acts on $\mP^1$
by the standard action of $SL_2(\mc)$ on $\mc^2$.
Then the irreducible components of
$\cz_X$ have the form $\fB\cdot (\z,\infty)\setminus (\z,\infty)$,
where $\z$ ranges over all of $X^\fT$ and
$\infty$ denotes the point $[0,1]\in \mP^1$.
The complete description of the  \ec\
of a \rv\ $X$  and principal subvariety $Y$ proved in \cite{BC} is given by
\begin{theorem} $($cf. \cite{BC}$)$ If $X$ is a \rv, then there exists a graded
$\mc$-algebra isomorphism $\rho_X:H_\fT^*(X)\to \mc[\cz_X]$.
Furthermore, if $Y$ is a \psv\ and $\cz_Y$ is the (reduced)
affine curve $\cz_X\cap (Y\times \mc)$, then there is also
a graded
$\mc$-algebra isomorphism $\rho_Y:H_\fT^*(Y)\to \mc[\cz_Y]$ making
the diagram
\begin{eqnarray}\label{CD1}
 \begin{CD}
             H_\fT^*(X) @>{\rho_X}>> \mc[\cz_X] \\
                {i_Y^*}  @VVV          {{\overline{i_Y^*}}} @VVV  \\
              H_\fT^*(Y)    @>{\rho_Y}>> \mc [\cz_Y],
\end{CD}
\end{eqnarray}
commutative, where the vertical maps are natural restrictions.
Moreover, the horizontal maps are $\mc[v]$-module maps under the
standard $\mc[v]$-module structure on $\tec(X)$ and $\tec(Y)$ and the
$\mc[v]$-module structure on $\mc[\cz_X]$ and $\mc[\cz_Y]$ induced by
the second projection.
\end{theorem}

One easily sees from the definitions
that if $v$ is an affine coordinate on
$\mc= \mP^1 \setminus [1,0]$, then $\mc[\cz_X]/(v)\mc[\cz_X] \cong A(Z)$.
Hence, Theorem 1 allows one to see the isomorphism
$H^*(X)\cong A(Z)$ in a natural way from
elementary properties of \ec.
Indeed, it turns out that
$\rho_X$ maps the augmentation ideal
$(v)H_\fT^*(X)$ in $H_\fT^*(X)$ to $(v)\mc[\cz_X]$,
so $A(Z)\cong H^*(X)$ since $H_\fT^*(X)/(v)H_\fT^*(X)\cong H^*(X)$.
Similarly, $H^*(Y)\cong \mc[\cz_Y]/(v)\mc[\cz_Y]$ if $Y$ is principal.
However, it is not in general
true that $\mc[\cz_Y]/(v)\mc[\cz_Y]$ is isomorphic to $A(Y\cap Z)$.
%See \cite{CRELLE}

Since surjectivity of $i_Y^*:H^*(X)\to H^*(Y)$ is equivalent to surjectivity of
$i_Y^*:H^*_\fT(X)\to H^*_\fT(Y)$,  (\ref{CD1})
suggests an approach to the surjectivity question.
Namely, since ${\overline{i_Y^*}}:\mc[\cz_X]\to \mc[\cz_Y]$
is surjective for any \bi\ subvariety $Y$,
the question of surjectivity of ${i_Y^*}$ boils down to
determining if there exists
an injective map $\rho_Y : H_\fT^*(Y)  \to \mc [\cz_Y]$
\st (\ref{CD1}) is commutative.
We will resolve this question by showing
\begin{theorem}\label{class-thm} A \bi\
subvariety $Y$ of a \rv\ $X$ is \pr\
\iff $H_\fT^*(Y)$ $($equivalently, $H^*(Y)$$)$
is generated by Chern classes of $\fB$-\eq\
algebraic vector bundles on $Y$.
%In particular, $\mc[\cz_Y] \cong H_\fT^*(Y)$  exactly when
%$H_\fT^*(Y)$ is generated by $\fB$-\eq\ algebraic vector bundles.
\end{theorem}

An example of a
\bisv\ which is not \pr\ is easily obtained. In fact, let $X$
be $\fB$-regular, and let $Y$ denote the union of all the
$\fB$-stable curves in $X$. Then $X^\fT = Y^\fT$, and it is
not hard to see that $i^*_Y$ is not surjective
if $\dim X>1$. However, $\tec(Y)$ is not generated
by $\fB$-\eq\ vector bundles. We will verify this
claim for $X=\mP^2$ in Section 5.
%Of course, $Y$ is neither normal nor irreducible.

%Let us begin with a basic example. The projective line ${\mP}^1$  under the action
%$$\begin{pmatrix}t&u\\ 0&t^{-1}\end{pmatrix}\cdot [x,y]
%= [ tx+uy, t^{-1}y]$$ is certainly a \rv.
Schubert varieties in $G/B$ or a $G/P$
are well known to be principal. Indeed, the
$B$-orbits are affine cells, and every Schubert
variety is a union of $B$-orbits.
Moreover, certain Springer fibres in $SL_n(\mc)/B$
are principal. To see this, first recall that a Springer fibre in $G/B$ is by
definition the fixed point set of a unipotent element of $G$.
Viewing $G/B$ as the variety of Borel subalgebras
of $\fg$, this definition is equivalent
defining a Springer fibre to be the set of all
Borels in $\fg$ containing a fixed nilpotent in $\fg$.
By classical result of Spaltenstein \cite{S}, the \coh\ map
$H^*(G/B)\to H^*(Y)$ is surjective for every Springer
fibre $Y$ provided $G=SL_n(\mc)$. Let $e$
denote a regular nilpotent in $\lie(B)$, and let
$\fB\circlearrowright G/B$
denote the regular action determined by $e$ as explained in
the first paragraph.
Any Springer fibre in $G/B$ corresponding
to a nilpotent in the centralizer $\fg^e$
of $e$ which is also a $\fT$-weight vector is $\fB$-invariant.
Thus such Springer fibres form a class of principal
subvarieties of $SL_n(\mc)/B$. This
example shows that principal subvarieties need not
be irreducible.

\section{Preliminaries}
Let $\la :\mc^* \to \fT$
and $\p:\mc \to \fU$ denote the one parameter subgroups
$$\la(t)=\begin{pmatrix}t&0\\ 0&t^{-1}\end{pmatrix}, \,\,\
\p(v)=\begin{pmatrix}1&v\\ 0&1\end{pmatrix}.$$
Then $\la(t)\p(v)\la(t)^{-1}=\p(t^2v)$ for all $t\in \mc^*$ and $v\in \mc$.
Suppose $X$ is a \rv\ with $X^\fU=\{o\}$, and note that $o\in X^\fT$. Put
$$X_o=\{x\in X\mid \lim_{t\to \infty}\la(t)\cdot x=o\}.$$
Recall from \cite{AC} (also see \cite{CRELLE})  that $X_o$ is an open
neighborhood of $o$ isomorphic with $T_o(X)$.
Thus, if $a_1, \dots, a_n$ ($n=\dim X$) are the weights
of $\la$ on $T_o(X)$ repeated with multiplicities, then all
$a_i<0$ and there exist affine coordinates
$u_1, \dots, u_n$ on $X_o$ which are quasi-homogeneous
of positive degrees $-a_1, \dots, -a_n$ with respect to $\la$.
That is,
$$\la(t)\cdot u_i=t^{-a_i}u_i,$$
for all $i$. The induced positive grading on  $\mc[X_o]=\mc[u_1, \dots, u_n]$
is frequently  called the {\em principal grading}.
Since the fixed point set $X^\fT$  of the
torus $\fT$ is finite  and contains $o$ (\cite{AC}),
we will write $X^\fT=\{\z_1, \dots, \z_r\},$  where $\z_1=o.$

%By \cite[Theorem 2]{AC}, each weight $a_i$ is even.
%Thus we have
%\begin{lemma}\label{weights-even} The $\fB$-action on a regular variety $X$ descends to
%an action of a Borel subgroup $\ob$ of $\mP SL_2(\mc)$.\end{lemma}
%\proof Since the weights $a_i$ are even, $\la(1)\cdot x=\la(-1)\cdot x$
%for all $x\in X_o$.
%Since $X_o$ is Zariski dense, it follows that  $\la(1)=\la(-1)$
%identically, so the action $\fB \circlearrowright X$
%descends to an action $\ob \circlearrowright X$, where
%$\ob =\fB/\{\pm I_2\}$. But $\ob$ is a Borel in $\mP SL_2(\mc)$. \qed

Note that the natural action of $SL_2(\mc)$ on $\mP^1$
induces a regular action. If $\fB$ denotes the upper
triangular matrices and $\fT$ the diagonals,
then $\z_1=[1,0]$, $\z_2=[0,1]$ and
the big cell is $\{[1,v]\mid v\in \mc\}$.

From now on, $\fB$ will denote the upper triangular Borel in $SL_2(\mc)$
and $\fT$ will be the diagonal torus; $X$ will always denote a \rv.
To simplify the notation, we will put $0=[1,0]$ and $\infty=[0,1]$
so $\mP^1=\mc\cup \{\infty\}$.
Notice that the diagonal action $\fB \circlearrowright (X\times \mP^1)$ is also
regular with $\fU$-fixed point $(o,0)$.
Define a projective curve $Z_X \subset X\times \mP^1$
as follows: let $Z_i$ be
Zariski closure of the orbit
$\fB (\z_i,\infty)$, where $i\ge 1$, and let
$$Z_X=\bigcup _{1\le i \le r} Z_i.$$
Thus $Z_X$ is a $\fB$-stable curve with $r=|X^\fT|$ irreducible components.
Moreover, the second projection
$p_2:X\times \mP^1 \to \mP^1$ induces an isomorphism
on each component. Hence $Z_X$ is
a bouquet of $r$ $\mP^1$s
passing through $(o,0)$. Now define an affine curve
\begin{equation}\label{fundcurve}
\cz_X =Z_X \cap (X_o \times \mc),
\end{equation}
where the intersection is assumed to be
in the sense of varieties, i.e. reduced.
This curve is $\fT$-stable but not $\fB$-stable.
In fact, $\cz_X$ is obtained from
$Z_X$ by removing the point at infinity
on each irreducible component.
In particular,  the coordinate ring $\mc[\cz_X]$
is a graded $\mc$-algebra via the principal grading, and
the projection $p_2$ induces a
$\mc[v]$-module structure on $\mc[\cz_X]$,
where $v$ denotes the affine coordinate on $\mc$.
For later use, we observe (cf. \cite[p.192]{BC})
\begin{equation}\label{box}
\makebox{if $v\ne 0$, then $(x,v)\in \cz_X$ \iff $\p(-v^{-1})\cdot x \in X^\fT$}.
\end{equation}

%\section{Cohomology of Regular Varieties}

We now recall the basic isomorphism $\rho_X:H_\fT^*(X)\to \mc[\cz_X].$
Since the odd \coh\ of a \rv\ is trivial,
the localization theorem in \ec\ implies
$i^*:H_\fT^*(X) \to H_\fT^*(X^\fT)$ is injective,
where $i:X^\fT \to X$ is the inclusion map. Thus, to each
$\A\in H_\fT^*(X)$, one may assign an $r$-tuple
$(\A_1, \dots, \A_r)\in {\bigoplus _i}\mc[v].$
Now define a function $\rho_X(\A)$ on $\cz_X$ by putting
\begin{equation}\label{defro}
\rho_X(\alpha) (x,v)=\alpha_{i}(v),
\end{equation}
if $(x,v)\in Z_i$. A key fact is that $\rho_X(\alpha)$ is a regular function on $\cz_X$.

%A key step in the proof is to show that
%$\rho_X(\A)$ is a regular function on $\cz_X$.
%This requires using the fact that $H_\fT^*(X)$
%is generated by the Chern classes of $\fT$-\eq\ vector bundles on $X$,
%which was observed in \cite{CRELLE} (see Lemma ???).

\section{Some applications of surjectivity}
We will now use Theorem 1 to draw some conclusions
about surjectivity. Suppose that $Y$ is  another \rv\ and
$F:Y\to X$ a $\fB$-equivariant map.
We claim $F$ induces a $\fT$-equivariant
map $\overline{F}:\cz_Y \to \cz_X$
by putting $\overline{F}(y,v)=(F(y),v)$. For $F(Y^\fT)\subset X^\fT$,
so if $(y,v)\in \cz_Y$, say $(y,v)=(b\cdot y_i, b\cdot \infty)$
for some $y_i \in Y^\fT$ and $b\in \fB$, we see that
$$(F(y),v)=(F(b\cdot y_i), b\cdot \infty)=(b\cdot F(y_i), b\cdot \infty)=
b\cdot (F(y_i), \infty)\in Z_X.$$
But $b\cdot (F(y_i), \infty)\ne (F(y_i), \infty)$
unless $b=1$, hence the claim. Thus there is a commutative diagram
\begin{eqnarray}\label{CD2}
\begin{CD} H_\fT^*(X) @>{\rho_X}>> \mc[\cz_X] \\
                {F^*}  @VVV          {\overline{F^*}} @VVV  \\
              H_\fT^*(Y)    @>{\rho_Y}>> \mc [\cz_Y].
\end{CD}
\end{eqnarray}

We now derive some consequences of this.
\begin{theorem}\label{thm1} Assume $F:Y\to X$ is a $\fB$-equivariant map
of \rvs\ varieties
\st the differential $dF_o$ of $F$ at $o$ is injective.
Suppose also that $F|Y^\fT$ is  injective.
Then the restriction map $F^*:H_\fT^*(X)\to H_\fT^*(Y)$
is surjective. Consequently, $F^*:H^*(X)\to H^*(Y)$
is also surjective.
\end{theorem}
\proof
First suppose $Y$ is a \bisv\ of $X$ which is also $\fB$-regular
and $F=i_Y$.
By Theorem 1, the morphism $\rho_Y$ in the diagram (\ref{CD2}) is an isomorphism.
Thus $i_Y^*$ is surjective. Now consider the general case.
Since $dF_o$ is injective,
it follows that $F(Y)$ is smooth at $o$. By the Borel Fixed
Point Theorem applied to the singular  locus of $F(Y)$,
it follows that $F(Y)$ is smooth, hence $\fB$-regular.
Hence $i_{F(Y)}^*: H_\fT^*(X)\to H_\fT^*(F(Y))$
is surjective. As $F$ is injective on $Y^\fT$ and $\fB$-\eq,
it follows that $\overline{F}: \cz_Y \to \cz_{F(Y)}$
is also injective. Since $F$ has a local holomorphic
inverse $G$ in a neighborhood of $o$, $G$ in fact
induces (by equivariant extension) a holomorphic
map $\overline{G}:\cz_{F(Y)} \to \cz_Y$ which is an
inverse (in the analytic category) to $\overline{F}$.
Note that  $F(Y)^\fT =F(Y^\fT)$; that is,
$F|Y^\fT$ is a bijection with $F(Y)^\fT$.
For, if $w\in F(Y)^\fT$,
then, by the Borel Fixed Point Theorem,
the subvariety $F^{-1}(w)$ of $Y$ contains a $\fT$-fixed point
due to the fact that it is $\fT$-invariant.
Therefore $\overline{F^*}: \mc[\cz_{F(Y)}] \to \mc[\cz_Y]$
is an isomorphism, giving the result. \qed

\medskip
Moreover, we also get
\begin{theorem} Two \pr\ subvarieties of a \rv\
with the same fixed point set have isomorphic \coh\ algebras
(both equivariant and classical). In particular, two such
subvarieties have the same dimension.
\end{theorem}
\proof This is an immediate consequence of Theorem 1. \qed

\medskip
Consequently, regular actions have a rather
special property.
\begin{corollary}\label{cor1} If $X$ is $\fB$-regular and $Y$ is a \bi\ subvariety
of $X$ such that $Y^\fT=X^\fT$, then either $Y=X$ or $Y$ is singular.
%In particular,
%there exist at most a finite number of
%strictly decreasing chains
%$$X\supset X_1\supset X_2 \supset \cdots \supset X_m,$$
%where each $X_i$ is a regular subvariety of $X$.
\end{corollary}
\proof By the previous theorem, if $Y$ is smooth, hence $\fB$-regular, then
$\dim Y=\dim X$. Since a \rv\ is necessarily irreducible,
$Y=X$. \qed

\medskip
Of course, there are examples of torus actions
on smooth projective varieties $Y\subsetneq X$ for which
$X^T=Y^T$ and the conclusion of Corollary \ref{cor1} fails.
 For example, if $X$ is the wonderful
compactification of a semi-simple algebraic group $G$ (over $\mc$)
with maximal torus $T$,
then all the fixed points of the (torus) action of $T\times T$
on $X$ lie on the unique closed $G\times G$-orbit $Y$
(cf. \cite{DP}). This gives
a proper smooth $T\times T$-stable subvariety $Y$ of $X$ for
which $X^{T\times T} =Y^{T\times T}$. It also shows that a
wonderful compactification  doesn't admit a regular action.
%For further details, see Section 4.

%An example of a singular \bi\ subvariety $Y$
%of $X$ where $Y^\fT=X^\fT$ is given by putting $Y=\fB X^\fT.$
%Thus, $Y$ is the union of all the $\fB$-stable curves in $X$.
%We will return to this example later.
%Since $X^\fB=\{o\}$, it does indeed follow that
%$Y=\fB X^\fT.$

Finally, we have
\begin{theorem} Let $X$ and $Y$ be \rvs\ or
\pr\ subvarieties, and let $F:Y \to X$ be a $\fB$-\eq\ map such that
$F(Y^\fT)= X^\fT$. Then $F^*: H_\fT^*(X) \to H_\fT^*(Y)$ is injective.
In particular, if  $F$ is surjective, then $F(Y^\fT)= X^\fT$,
so $F^*$ is injective on \ec.
\end{theorem}

\proof Since
$F(Y^\fT)=X^\fT$, $\overline{F}: \cz_Y \to \cz_X$ is also surjective.
Since $\cz_X$ and $\cz_Y$ are affine, it follows that
the comorphism $\overline{F}^*:\mc[\cz_{X}]\to \mc[\cz_{Y}]$
is injective. Therefore, $F^*: H_\fT^*(X) \to H_\fT^*(Y)$ is injective.
If $F$ is surjective, then the Borel Fixed Point Theorem implies
$F(Y^\fT)=X^\fT$.
\qed

\medskip
\ni
{\bf Remark.} If $X$ and $Y$ are smooth projective varieties
and $F:Y\to X$ is surjective, it is well known
that $F^*$ is always  injective on ordinary \coh.

\section{A remark on a formula of Kostant and macdonald}
Let $G$ be a complex semi-simple algebraic group, $B$ a Borel subgroup of $G$
and put $\fb=\lie(B)$.
We will now describe an interesting class of
$\fB$-regular subvarieties of $G/B$ which includes
all smooth Schubert varieties.
Let $e\in \fb$ be a regular nilpotent in $\fg$, and
recall from Section 1 that $e$ determines a regular action
$\fB\circlearrowright G/B$ on the flag variety of $G$
by left translation
such that the identity coset $B\in G/B$ is the
unique $\fU$-fixed point. To simplify notation, let us identify
$\fB$ with its homomorphic image in $G$
such that $e\in \lie(\fB)$
(see the comment in the first paragraph of the Introduction).
Let $\fh$ denote a $\lie(\fB)$-submodule of $\fg$
containing $\fb$, and put $Y_\fh=\overline{\exp(\fh)B}.$

\begin{lemma} $Y_\fh$ is a $\fB$-invariant subvariety of $G/B$.
\end{lemma}
\proof
Since the exponential map $\exp:\fg \to G$ is $G$-\eq\
for the adjoint action of $G$ on $\fg$, for any $x\in \lie(\fB)$, we have
$$ \exp \big(Ad\big(\exp (tx)\big) y\big)B=\exp(tx)\exp(y)\exp(-tx)B=\exp(tx)\exp(y)B.
$$
But for any $x,y\in \fg$,
$Ad\big(\exp (tx)\big)(y)=e^{ad(tx)} (y)$,
so the term on the left side of the above identity
is in $\exp(\fh)B$
if $y\in \fh$. Thus $Y_\fh$ is stable under $\fB$.\qed

\medskip
\noindent
{\bf Remark:} If $\fh$ is $B$-stable, then
$Y_\fh$ is also $B$-stable and hence is a Schubert variety.

\medskip
When the variety $Y_\fh$ is smooth, it is a $\fB$-regular variety and
one can now apply the formula in \cite{AC}
for the Poincar\'e polynomial of
a \rv\ $X$ to find
an interesting class of polynomials associated to any
complex semisimple Lie algebra $\fg$. Let us first recall the formula.
Let $a_1, \dots, a_n$ ($n=\dim X$) denote the weights
of $\la$ on $T_o(X)$ introduced in Section 2,
and recall the $a_i$ are negative integers.
Then
\begin{equation}\label{prodform}P(X, t^{1/2})=\prod_{1\le i\le n}
\frac{(1-t^{-a_i+1})}{(1-t^{-a_i})}.
\end{equation}
In the case  $X=Y_\fg=G/B$,
the weights $a_i$ are the negatives of the heights of the
positive roots with respect
to any maximal torus $T$ of $G$ contained in $B$, so
we recover a well known formula
$$P(G/B,t^{1/2})=\prod_{\A>0}\frac{(1-t^{ht(\A)+1})}{(1-t^{ht(\A)})}$$
of Kostant and Macdonald.

\section{A lemma and an example}

We will now prove some facts about \bi\ subvarieties
of a \rv\ $X$. First note

\begin{lemma}\label{psi-lemma} Let $Y$ be  a \bi\ subvariety of $X$
with vanishing odd \coh, and let
$I(\cz_Y)\subset \mc[\cz_X]$ denote the ideal of $\cz_Y$. Then:

\medskip
\begin{itemize}
\item[(i)] $\rho_X(H_\fT^{ev}(X,Y))=I(\cz_Y)$; and

\medskip
\item[(ii)]
there exists a $\mc[v]$-algebra isomorphism
$$\psi_Y: \mc[\cz_Y] \to i_Y^*\big(H_\fT^*(X)\big) \subset H_\fT^*(Y).$$
\end{itemize}
In fact, $\psi_Y=i_Y^*\rho_X^{-1} (\overline{i_Y^*})^{-1}.$
\end{lemma}

\proof
We first show that putting  $\psi_Y=i_Y^*\rho_X^{-1} (\overline{i_Y^*})^{-1}$
gives a well defined map.
Since $\cz_Y$ is the union of the irreducible
components of $\cz_X$ which meet $Y\times \mc,$
the indeterminacy introduced
by $(\overline{i_Y^*})^{-1}$ is supported on the complement of $\cz_Y$.
Thus $\psi_Y$ is indeed well defined since
$\overline{i_Y^*}$ is surjective,
and two classes in $H_\fT^*(X)$ which have the same
image under $\overline{i_Y^*}\rho_X$ have the same image under
$i_Y^*$, by the localization theorem and the definition of $\rho_X$.

We now show $\rho_X\big(H^{ev}_\fT(X,Y)\big)=I(\cz_Y)$.
Since $\cz_Y$
is an affine curve with $|Y^\fT|$ components, and $p_2:\cz_Y\to \mc$
is a flat map whose restriction to each component is an isomorphism
(by \cite[Prop. 2]{BC}), the rank of $\mc[\cz_Y]$ over $\mc[v]$ is $|Y^\fT|$.
Furthermore, by the long exact sequence
of \coh\ and the localization
theorem, the rank of $H^{ev}_\fT(X,Y)$ is
$|X^\fT|-|Y^\fT|$. But the rank
of $\ker \overline{i_Y^*}$ is also
$|X^\fT|-|Y^\fT|$. As shown above,
$\ker \overline{i_Y^*}\subset \rho_X(H^{ev}_\fT(X,Y))$,
so it follows that
$\overline{i_Y^*}( \rho_X(H^{ev}_\fT(X,Y)))$ has rank zero
in the free module $\mc[\cz_Y],$ hence is trivial.
Therefore, $I(\cz_Y)= \rho_X(H^{ev}_\fT(X,Y))$.
To finish, we only need to show $\psi_Y$ is injective.
But this follows immediately from part (i).
\qed

\medskip
The following example shows that non-\pr \
subvarieties exist. Let
$$\la(t)=\begin{pmatrix}t^2 &0&0\\0&1&0\\0&0&t^{-2}\end{pmatrix}$$
and
$$\p(v)=\begin{pmatrix}1&v&v^2/2\\ 0&1&v\\0&0&1
\end{pmatrix}.$$
Then $\la(t)\p(v)(\la(t))^{-1}=\p(t^2v)$, so $\la$
and $\p$ determine a  two dimensional solvable
subgroup of $SL_3(\mc)$, and hence a
regular action
$\fB \circlearrowright \mP^2$.
The $\fT$ fixed points are $o=[1,0,0], ~\z_2=[0,1,0],$
and $\z_3=[0,0,1]$. Let $w_1$ and $w_2$
be the usual affine coordinates around $o$.
Then the closures of $Y_1=\{w_2=0\}$ and $Y_2=\{2w_2=w_1^2\}$ are the
two $\fB$-curves in $\mP^2$. Let $Y=\overline{Y_1}\cup \overline{Y_2}$.
Then $Y$ is a
\bisv\ such that $H^*(\mP^2)\to H^*(Y)$ is not surjective.
Indeed, $\dim H^2(\mP^2)<\dim H^2(Y)$. On the other hand,
$H_\fT^*(Y)$ is not generated by Chern classes of equivariant vector bundles.
To see this, note that $Y$ has vanishing odd \coh,
so one can use the localization theorem to compute its
\ec. In fact, by a well known result of Goresky, Kottwitz and MacPherson
\cite{GKM} (also see \cite{BRION}), the image of $H_\fT^*(Y)$
in  $H_\fT^*(Y^\fT)$ consists of triples $(f_1(t),f_2(t),f_3(t))$
with all $f_i\in \mc[t]$ \st $f_1(0)=f_2(0)=f_3(0)$.
But not all classes of this form
arise as polynomials in Chern classes of \eq\ line bundles
on $Y$, due to the fact that equivariance forces the further
condition  $f_2=f_3=-f_1$.

\medskip
\noindent
{\bf Remark.} We do not know of an example of an irreducible
\bisv\ of a \rv\ which is not principal.

\section{Equivariant Chern classes}

Let $Y$ be a \bisv\ of a \rv\ $X$.
The purpose of this section is to
consider when the fundamental isomorphism
 $\rho_X:\tec(X) \to \mc[\cz_X]$ defined in (\ref{defro})
can be defined for $Y$. For this, we need to consider
$\fB$-\eq\ vector bundles on $Y$.

Recall that if $V$ is an algebraic variety with an action of an algebraic group $G$,
then a $G$-linearization of an algebraic vector bundle $E$ on $V$
is an action $G\times E\to E$ ($(g,h)\to g\cdot h$) \st
for all $y\in V$, the restriction of
$g\in G$  is a $\mc$-linear map $E_y\to E_{g \cdot y}$.
In particular if $y \in V^G$, one obtains a representation of $G$ on $E_y$
and hence a representation of $\lie(G)$ on $E_y$. For $\xi \in \lie(G)$, we let
$\xi_y$ denote the corresponding endomorphism of $E_y$.

Suppose the vector bundle $E$ admits a $G$-linearization.
Recall that the $k$-th
equivariant Chern class $c_k^G(E)\in H^{2k}_G(V)$ of $E$
is the $k$-th Chern class of the vector bundle $E_G=
(E\times {\mathcal E})/G$, where ${\mathcal E}$ is a contractible
free $G$-space. The restriction of $c_k^{G}(E)$ to each $y \in V^G$ is
the polynomial on $\lie(G)$ defined as
\begin{equation} \label{eqn-Chern-class}
c_k^G(E)_y(\xi) = \textup{Tr}_{\bigwedge^k E_y}(\xi_y).
\end{equation}

We now turn to the case  $G=\fB$, $X$ $\fB$-regular and $Y$ a \bisv.
Put
$$ \cw = \begin{pmatrix}1&0\\ 0&-1\end{pmatrix}, \,\,\
\cv =\begin{pmatrix}0&1\\ 0&0\end{pmatrix}.$$
By \cite[Lemma 1]{BC}, if $E$ is $\fB$-\eq\ and $(y,v)\in \cz_X$, then
$$\rho_X(c_k^{\fB}(E))(y,v)=\textup{Tr}_{\bigwedge^k E_y}(v\cw -2\cv)_y,$$
where  $(v\cw -2\cv)_y$ is the endomorphism
determined by $v\cw -2\cv\in \lie(\fB).$
This is a key fact since the smoothness of $X$ implies $\tec(X)$
is generated by Chern classes of $\fB$-\eq\ vector bundles on $X$
\cite[Prop. 3]{CRELLE}.

Let $\ch (Y)$ denote the subalgebra
of $\tec(Y)$ generated by \eq\
Chern classes of $\fB$-equivariant vector bundles on $Y$.
As just noted, $\ch(X)=\tec(X)$.
The goal of the remainder of this section is to prove that if
$E$ is a $\fB$-linearized
vector bundle on $Y$ then $\rho_Y(c_k^\fB(E))$
is a regular function on the affine
curve $\cz_Y$. Thus, there exists a well defined map
$$\rho_Y:\ch(Y) \to \mc[\cz_Y]$$
\st the diagram
\begin{eqnarray}\label{CD3}
\begin{CD} H_\fT^*(X) @>{\rho_X}>> \mc[\cz_X] \\
                {i_Y^*}  @VVV          {\overline{i_Y^*}} @VVV  \\
              \ch(Y)    @>{\rho_Y}>> \mc [\cz_Y].
\end{CD}
\end{eqnarray}
is commutative.

Before beginning the proof, assume $v\ne 0$ and define
$$\fH_v=\p(1/v) \fT \p(-1/v).$$
By (\ref{box}), $(y,v)\in \cz_Y$ \iff $y\in Y^{\fH_v}$.
Also, put
$$\mathcal{H} = \overline{\bigcup_{v \neq 0} \fH_v \times \{v\}}
\subset \fB \times \mc.$$
%Recall $\fB\subset \mP SL_2(\mc)$.

\begin{lemma} Let $\pi:\mathcal{H} \to \mc$ be projection on the second factor. Then:
$$\pi^{-1}(0) = \fU\cup -(\fU).$$
\end{lemma}
\proof
First note the identity
\begin{equation}\label{mat-id}
\begin{pmatrix} 1 & v^{-1} \\ 0 & 1 \end{pmatrix} \begin{pmatrix} a & 0\\0 &
a^{-1}
\end{pmatrix}
\begin{pmatrix} 1 & -v^{-1} \\ 0 & 1 \end{pmatrix} = \begin{pmatrix} a & (1-a^2)(av)^{-1} \\ 0 &
a^{-1}
\end{pmatrix}.
\end{equation}
One must find all the possible finite
limits of the right-hand side as $v \to 0$.
Take a sequence $(a_i, v_i)\ne (0,0)$  such that
$v_i \to 0$, and the
right-hand side of (\ref{mat-id}) has a finite limit.
In order for this to happen, $(1-a_i^2) / a_i \to 0$,
so $a_i^2 \to 1.$ Hence, if a limit exists, it lies in
$\fU\cup -(\fU)$. The reverse inclusion is similar.
\qed

%\medskip
%\ni {\bf Remark.} If we replace the action of $\fB$ by the induced
%action of the associated Borel $\ob$ in $\mP SL_2(\mc)$, then the fibre
%$\pi(0)$ is just $\fU$.

\medskip
Consequently, if we put  $\fh_v = \lie(\fH_v)$ for $v\ne 0$ and set
$\fh_0 = \lie(\fU)$, then the family of Lie algebras of the fibres of ${\mathcal H}$ is
$$\fh =\bigcup_{v\in \mc} \fh_v \times \{v\} \subset \lie(\fB)\times \mc,$$
and, moreover, $v\cw-2\cv$ is a non-vanishing regular section
of $\pi$. For convenience, let $s(v)=v\cw-2\cv$.

Now let $E$ be a $\fB$-linearized vector bundle on $Y$,
and let $\tilde{E}$ be the pull-back of $E$ to $Y \times \mc$.
Assume $v \in \mc$ and $y \in Y^{\fH_v}$. That is, $(y,v)\in \cz_Y$.Then
$E_y$ admits an $\fH_v$-representation and hence an $\fh_v$-representation.
For $\xi \in \fh_v$ let $\xi_y$ denote the corresponding endomorphism on $E_y$.
Then the upshot of our discussion is the following:
\begin{lemma} \label{lemma-Tr-regular}
The map from $\cz_Y$ to $\End(\tilde{E}_{|\cz_Y})$ given by
$$(y, v) \mapsto s(v)_y,$$ is a regular section
of the bundle $\End(\tilde{E}_{|\cz_Y})$.
Consequently the map
$$(y, v) \mapsto \Tr_{\bigwedge^k E_y}(s(v)_y)$$ is a
regular function on $\cz_Y$.
\end{lemma}

\begin{lemma} \label{lemma-Chern-class-regular}
Let $E$ be a $\fB$-linearized vector bundle on $Y.$ Then
$\rho_Y(c_k^\fT (E))$ is a regular function on $\mc[\cz_Y]$.
\end{lemma}
\begin{proof}
View $v$ as an element of $\lie(\fT)$. As above let $v_{\zeta_j}$ denote the
corresponding endomorphism on $E_{\zeta_j}$.
By Equation (\ref{eqn-Chern-class}) we have
$$c_k^{\fT}(E)_{\zeta_j}(v) = \Tr_{\bigwedge^k E_y}(v_{\zeta_j}).$$
Put $\zeta_j(v) = \p(v^{-1})\zeta_j$. Then the $\zeta_j(v)$ are the fixed points
of $\fH_v$. Recall that $s(v) \in \fh_v$, hence it gives an endomorphism
$s(v)_{\zeta_j(v)}$ of the fibre $E_{\zeta_j(v)}$. Since $E$ is $\fB$-linearized, the
element $\p(v^{-1})$ gives an isomorphism $E_{\zeta_j}\to E_{\zeta_j(v)}$; moreover
the endomorphism $v_{\zeta_j(v)}$ on $E_{\zeta_j(v)}$ is conjugate to the endomorphism
$v_{\zeta_j}$ on $E_{\zeta_j}$. Put $y = \zeta_j(v)$. Then $(y, v)\in \cz_Y$
and
$$\Tr_{\bigwedge^k E_{\zeta_j}}(v_{\zeta_j}) = \Tr_{\bigwedge^k E_{y}}(s(v)_y).$$
But by Lemma \ref{lemma-Tr-regular}, the function $(y, v) \mapsto \Tr_{\bigwedge^k E_y}(s(v)_y)$
is a regular function on $\cz_Y$.
\end{proof}

To summarize the discussion of this section,
we state
\begin{theorem} If $Y$ is a \bisv\ of a
\rv\ $X$, then we obtain a $\mc[v]$-algebra homomorphism
$\rho_Y:\ch(Y) \to \mc[\cz_Y]$ \st the diagram
(\ref{CD3}) is commutative.
\end{theorem}

\section{Classification of \pr\ subvarieties}
We now classify \pr\ subvarieties. First, we describe the image of  $\tec(X)$
in $\tec(Y)$ for a \bisv\ $Y$ with vanishing odd \coh.

\begin{theorem}\label{last-thm} Suppose $Y$ is \bisv\ of a \rv\ $X$
with vanishing odd \coh. Then

\medskip
\begin{itemize}
\item[(i)] $\rho_Y:\ch(Y) \to \mc[\cz_Y]$ is an isomorphism, and \\

\item[(ii)]  $i_Y^*(\tec(X))=\ch(Y)$. \\
\end{itemize}
Consequently,  $\ch(Y)$  is exactly
the subalgebra generated by Chern classes of $\fB$-\eq\ vector bundles on $Y$
which are pull backs of $\fB$-\eq\ vector bundles on $X$.
\end{theorem}
\proof  It is clear from the commutativity of
(\ref{CD3}) that $\rho_Y$ is surjective.
The definition of $\rho_Y$ and localization
(applied to $Y$) implies it is injective.
This proves (i). We now show (ii).
By Lemma \ref{psi-lemma}, we have an isomorphism
$\psi_Y:\mc[\cz_Y]\to i_Y^*(\tec(X))$.
We claim that $\rho_Y \psi_Y=1$.
In fact, by (\ref{CD3}) and the definition of $\psi_Y$,
$$\rho_Y \psi_Y=\rho_Y i_Y^*  \rho_X^{-1}\overline{i_Y^*}^{-1}
= \overline{i_Y^*}\rho_X \rho_X^{-1} \overline{i_Y^*}^{-1},$$
which is clearly the identity. Thus
$$\rho_Y (i_Y^*(\tec(X)) =\mc[\cz_Y],$$
which certainly implies (ii).
\qed

\medskip
\begin{corollary} A \bisv\ $Y$ of a \rv\ $X$ is \pr \ \iff $\ch(Y)=\tec(Y)$.
\end{corollary}
\proof The necessity is clear. But if $\ch(Y)=\tec(Y)$, then
$Y$ has vanishing odd \coh, so the result follows from
Theorem \ref{last-thm} (b). \qed

\medskip
Theorem \ref{class-thm} follows immediately from this corollary.
We conclude with an application.
\begin{corollary}
Suppose $Y$ is a normal \bisv\ of a \rv\ $X$ \st $H^*(Y)$ is generated by Chern classes of line
bundles. Then $Y$ is principal.
\end{corollary}

\proof Let $L$ be an algebraic line bundle on $Y$. By \cite{KKL}, some power $L^m$
is $\fB$-\eq. Hence the first Chern classes of $\fB$-\eq\ lines bundles on $Y$
generate $H^*(Y)$. This is equivalent to saying that $\ch(Y)=\tec(Y)$. \qed

\bigskip
\ni
{\bf Acknowledgement:} We would like to thank Michel Brion and Jochen Kuttler
who each noticed an error in a previous version of this paper.

\newpage
{\footnotesize

\bigskip
\noindent James B. Carrell, University of British Columbia,
    Vancouver, B.C., Canada \\{\it Email address:} {\sf
carrell@math.ubc.ca}\\

\noindent Kiumars Kaveh, University of Toronto, Toronto,
    ON., Canada \\{\it Email address:} {\sf kaveh@math.toronto.edu}\\

\end{document}